\documentclass[12pt, twoside, leqno]{article}
\usepackage{latexsym}
\usepackage{amsmath}
\usepackage{amsfonts}
\usepackage{amssymb}
\usepackage{amsthm}
\usepackage{color}
\numberwithin{equation}{section}

\input xy
\xyoption{all} \hfuzz=1.5pt \tolerance=500
\theoremstyle{plain}

\theoremstyle{definition}

\newcommand{\tr}{\triangleright}
\newcommand{\tl}{\triangleleft}
\newcommand{\cd}{{\cdot}}
\newcommand{\ot}{{\otimes}}

\newcommand{\la}{\langle}
\newcommand{\ra}{\rangle}
\newcommand{\id}{{\bf 1}}
\newcommand{\di}{\diamondsuit}
\newcommand{\q}{\quad}

\newcommand{\va}{\varphi}
\newcommand{\rr}{{\cal R}}
\newcommand{\ii}{\infty}
\newcommand{\mt}{\mapsto}
\newcommand{\co}{{\mathbb C}}

\newcommand{\bb}{{\cal B}}

     \newcommand{\lm}{{\lambda}}

\newcommand{\N}{{\mathbb N}}  
  
\def\x{\relax\ifmmode {\mbox{*}}\else*\fi}

\newcommand{\ed}{\end{document}}

\begin{document}


\bigskip
\centerline{{\bf{\large The existence of $p$-convex tensor products of }}}

\centerline{{\bf{\large $L_p(X)$--spaces for the case of an arbitrary measure }}}\footnote{This 
paper was written with the support of the Russian Foundation for Basic Research (grant no. 
19-01-00447).}

\bigskip
\centerline{A.~Ya.~Helemskii  }

\bigskip
\begin{abstract}

\medskip 
We obtain a far-reaching generalization (in several directions) of the theorem of A.~Lambert on the 
existence of the projective tensor product of operator sequence spaces. This result is obtained in 
the context of spaces, generalizing $p$-multinormed spaces of Dales et al. which are based on an 
arbitrary, perhaps non-discrete measure. 

\end{abstract}



\section{Introduction}

The main result of the present paper is the following theorem. All notions that participate in it, 
will be gradually explained. Below $X$ denotes an arbitrary measure space which is not atomic with 
finite set of atoms and which is supposed, for simplicity, to be separable. 

\medskip
{\bf The main Theorem.} {\it An arbitrary pair of near-$L_p(X)$-spaces, where $1\le p<\ii$, has a 
$p$-convex tensor product.}

\medskip
A far-away predecessor of this result is a theorem of A.~Lambert~\cite[\S 3.1.1]{lam} on projective 
tensor products of his ``operator sequence spaces''; the latter are situated, in a sense, between 
normed spaces and abstract operator spaces. Afterwards, a group of mathematicians  (Dales, 
Polyakov, Daws, Pham, Ramsden, Laustsen, Oikhberg, Troitsky; see~\cite{dal} and 
also~\cite{dpol,ddaws, ddaws2}) introduced more general structures than Lambert spaces, called 
$p$-multinormed spaces; $p\in[1,\ii]$. After their papers, in the frame-work of the so-called 
non-coordinate approach, the $L$-spaces were introduced in~\cite{hel}. The latter, provided 
$L:=L_p(X)$, can be considered as ``multinormed spaces based on arbitrary measures''. Indeed, in 
the case of $X:={\Bbb N}$ with the counting measure they transform to $p$-multinormed spaces 
of~\cite{dal} (and if, in addition, $p:=2$, to Lambert spaces). 

``The main theorem'' was proved in~\cite{hel} under the additional condition that the set of atoms 
in $X$ is either empty or infinite. For a time there was a suspicion that for an arbitrary $X$ the 
theorem is false. (It was based on the known bad properties of $L_p(X)$ as a tensor factor, see 
~\cite[\S12.1]{def}). However, recently T.~Oikhberg kindly sent to the author a 
preprint~\cite{oik}, where he has constructed an isomorphism between the categories of 
$\ell_p$-spaces and $L_p(X)$-spaces for arbitrary $X$. Apart from the independent value of this 
result, its proof was based on a construction that, as it happened, allowed us to dispense of the 
afore-mentioned additional condition on $X$. As a matter of fact, one can even consider tensor 
products of the so-called ``stratified spaces'', more general than $L_p(X)$-spaces. But this is 
outside of the scope of this paper. 

The author is indebted to T.~Oikhberg and N.~T.~Nemesh for valuable discussions.

\section{$L$--spaces and $L$--boundedness. }

As usual, we denote by $\bb(E)$ the space of all bounded operators on a normed space $E$, endowed 
with the operator norm. 
Two projections $P$ and $Q$ on $E$ are called {\it transversal}, if $PQ=QP=0$. The symbol $\ot$ is 
used for the algebraic tensor product of linear spaces and linear operators, and also for  
elementary tensors. 
 
Choose and fix (so far arbitrary) normed space $L$, which we shall call the {\it bаse space}. 

Our principal example  of a base space is $L_p(X)$, where $1\le p<\ii$, and $X$ is a measure space, 
which is not reduced to a finite set of atoms, or, equivalently, $L_p(X)$ is infinite-dimensional. 
To make our text shorter, we shall always assume that all our measures have a countable basis. 

The {\it amplification} of a given linear space $E$ is the tensor product $L\ot E$. Usually we 
shortly denote it by $LE$, and an elementary tensor, say $\xi\ot x; \xi\in L, x\in E$, by $\xi x$. 
Note that $LE$ is a left module over the algebra $\bb$ with the outer multiplication `` $\cd$ '', 
well defined by $a\cd(\xi x):=a(\xi)x$.

\medskip
{\bf Definition 2.1.} A norm on $LE$ is called {\it $L$--norm on} $E$, if the left $\bb(L)$-module 
$LE$ is contractive, that is if we always have the estimate $\|a\cd u\|\le\|a\|\|u\|$. This 
estimate, as well as an equivalent estimate $\|a\ot\id_E\|\le\|a\|$, will be called {\it 
contractibility property.} The space $E$, endowed by an $L$--norm, is called {\it ${\bf 
L}$--space}. If we only know that we have the indicated estimate  for operators of rank 1, we speak 
about near-$L$-norms and, accordingly, about near-$L$-spaces. 

\medskip
{\bf Remark 2.2.} The class of near-$L$-spaces that we shall need in the proof of the main theorem, 
is bigger than the class of $L$-spaces. Let $L:=L_p(X)$, $E$ be a normed space, and the norm on 
$LE$ is given by the identification of that space with the corresponding dense subspace in 
$L_p(X,E)$. Then we obtain a near-$L$-space but, generally speaking, not an $L$-space: the 
indicated estimate fails already in the case $p=2, E=\ell_1$ (see~\cite[p.147]{def}). 

\medskip
 As to the papers, cited above, they consider, after the translation into the 
``index-free'' language, the case $X:=\N$ with the counting measure. In particular, the notion of 
an $L_p(X)$-space for such an $X$ is equivalent to that of a $p$-multinormed space in~\cite{dal}. 

\medskip
{\bf Proposition 2.3.} {\it Suppose that $E$ is a normed space, and a cross-norm $\|\cd\|$ is given 
on $LE$. Then $\|\cd\|$ is a near-$L$-norm on $E$ iff for all $f\in L^*$  we have 
$\|f\ot\id_E\|\le\|f\|$.} 

\smallskip
$\tl$ Every $a\in\bb(L)$ of rank 1 acts as $\eta\mt f(\zeta)\xi$ for some $\xi$ and $f$, so $\|a 
\|=\|f\|\|\xi\|$. It is easy to verify that $a\cd u=\xi\Big((f\ot\id_E)u\Big)$, for all $u\in LE$. 
Therefore the estimates $\|a\cd u\|\le\|a\|\|u\|$ and $\|(f\ot\id_E)u\|\le\|f\|\|u\|$ are 
equivalent. $\tr$ 

\medskip
A near-${L}$--space $E$ becomes a normed space in the ``classical'' sense, if  for $x\in E$, we set 
$\|x\|:=\|\xi x\|$, where $\xi\in L$ is an arbitrary vector with $\|\xi\|=1$. Clearly, the result 
does not depend on the choice of $\xi$. The obtained normed space is called the {\it underlying 
space} of a given ${L}$-space, and the latter is called an {\it ${L}$--quantization} of a former. 
We use such a term by analogy with quantizations in operator space theory; see, 
e.g.,~\cite{ef5},~\cite{efr} or~\cite{heb2}. 

\medskip
It is easy to verify that the complex plane $\co$ has the only $L$--quantization, given by the 
identification of $L\co$ with $L$. However, as a rule, general normed spaces have a lot of 
${L}$--quantizations. In particular, by endowing $LE$ with the norm of (non-completed) projective, 
respectively injective tensor product of normed spaces, we obtain two, generally speaking, 
different $L$-quantization, called maximal, respectively minimal. (See~\cite{hel} for details.) 

Suppose we are given an operator $\va:E\to F$ between linear spaces. Denote, for brevity, the 
operator $\id_L\ot\va:LE\to LF$ (taking $\xi x$ to $\xi\va(x)$) by $\va_\ii$ and call it {\it 
amplification} of $\va$. Obviously, $\va_\ii$ is a morphism of left $\bb(L)$-modules. 

\medskip
{\bf Definition 2.4.} An operator $\va:E\to F$ between $L$--spaces is called {\it ${\bf 
L}$--bounded} or {\it $L$--contractive}, if the operator $\va_\ii$ is bounded or contractive, 
respectively. 

As to numerous examples and counterexamples see, e.g.,~\cite{hel}, and also~\cite{dal}.

\medskip 
To define amplifications of bilinear operators, we need a certain additional structure, called in 
what follows {\it $\di$--operation} or ``diamond operation'' on $L$. This is a bilinear operator 
$\di:{L}\times L\to L$ of norm one. We shall write $\xi\di\eta$ instead of $\di(\xi,\eta)$. 

For ``most'' $X$, $L_p(X)$ has a natural, in a sense, diamond operation (see~\cite[\S3]{hel}, and 
also, in the case of a discrete measure,~\cite[\S 1.2.2]{lam}). But we emphasize that our main 
theorem is valid for arbitrary $\di$.  

\medskip
Now let $\rr:E\times F\to G$ be a bilinear operator between linear spaces. Its {\it amplification} 
is the bilinear operator $\rr_\ii:LE\times LF\to LG$, well-defined (because of the bilinearity) on 
elementary tensors by $\rr_\ii(\xi x,\eta y)=(\xi\di\eta)\rr(x,y)$ . 

\medskip
{\bf Definition 2.5}. A bilinear operator $\rr$ between ${L}$--spaces is called {\it ${
L}$--bounded} or {\it ${L}$--contractive},  if its amplification is (just) bounded, or 
contractive, respectively.

\medskip
In the case $L=\ell_2$ and a particular $\di$, taking sequences $\{\xi_n\}$ and $\{\eta_n\}$ into 
(arbitrarily renumerated) double sequence $\{\xi_n\eta_m\}; m,n\in{\Bbb N}$, we obtain, in 
equivalent terms, the definition of an $L$-bounded bilinear operator, given by Lambert. 

Again, see~\cite{hel} for numerous examples.

 \section{$p$-convex tensor product and preliminaries of its existence}

\medskip
{\bf From now on, and up to the end of the paper, we assume that $L:=L_p(X); p\in[1,\ii)$} (i.e. we 
are within the context of our main example of $L$), {\bf and that we fix an arbitrary 
$\di$-operation on our base space.} 

Let $Y$ be a measurable subset in $X$. Consider the projection $P_Y\in\bb(L)$, acting as $f\mt 
f\chi$, where $\chi$ is a characteristic function of $Y$. A projection of that kind will be called 
{\it proper}. Clearly, two proper projections are transversal iff the intersection of the 
respective measurable subsets has measure 0. 

\medskip
Let $E$ be a linear space. We call a projection $P\in\bb(L)$ a {\it support } of an element $u\in 
LE$, if $P\cd u=u$. 

\medskip
{\bf Definition 3.1.} A near-$L$--space $E$ is called {\it $p$--convex,} if for any $u,v\in LE$, 
with transversal proper supports, we have $\|u+v\|\le(\|u\|^p+\|v\|^p)^{\frac{1}{p}}$. 

\medskip
The introduced class of $L$-spaces, being a generalization for arbitrary $p$ of column operator 
spaces, is, in our opinion, the most interesting. For  the special case $L:=\ell_p$, the given 
definition is equivalent to  the definition of a $p$--convex $p$--multinormed space, given 
in~\cite{dal}. Also it worth mentioning, in this connection, the theory of $p$--operator spaces of 
Daws~\cite{daw}; see also earlier papers of Pisier~\cite{pi2} and Le Merdy~\cite{lem}.  

As an example, one can easily show that every $L$-space with the minimal quantization is 
$p$-convex. Another example is provided by the near-$L$-space from Remark 2.2.

Now let $E$ and $F$ be two {\it arbitrary chosen} near-$L$--spaces.

\medskip
{\bf Definition 3.2.} A pair $(\Theta,\theta)$, consisting of a $p$-convex $L$-space $\Theta$ and 
an $L$--contractive bilinear operator $\theta:E\times F\to\Theta$, is called (non-completed) {\it 
$p$-convex tensor product of $E$ and $F$} if, for every $p$-convex $L$-space $G$ and every 
${L}$--bounded bilinear operator $\rr:E\times F\to G$, there exists a unique ${L}$--bounded 
operator $R:\Theta\to G$ such that the diagram 
\[
\xymatrix@R-10pt@C+15pt{
E\times F \ar[d]^{\theta} \ar[dr]^{\rr} & \\
\Theta \ar[r]^R  &  G  } 
\]
\noindent is commutative, and moreover we have $\|R_\ii\|=\|\rr_\ii\|$. 

In what follows, the property of the pair in question will be called {\it the universal property}.

\medskip
We emphasize that $\Theta$ and $G$ are supposed (in comparison to $E$ and $F$) to be $L$-spaces, 
and not just near-$L$-spaces. 

\medskip
{\bf Remark 3.3}. We see that the $L$--spaces $\Theta$ and $G$ are assumed to be $p$--convex. Other 
assumptions lead to other types of tensor products. For instance, if we shall take the class of all 
$L$-spaces, we shall come to an essentially different concept, the so-called general tensor product 
of our  $E$ and $F$. This variety has its own existence theorem; this is Theorem 4.6 in~\cite{hel}. 
Nevertheless $p$-convex tensor products, being in the case $p=2$ intimately connected with the 
projective tensor products of operator spaces, discovered by Blecher/Paulsen~\cite{blp} and 
Effros/Ruan~\cite{er2}, seem to be most interesting. 

\bigskip 
Thus, all notions that participate in the formulation of our main theorem, are explained, and we 
can proceed to its proof.

As it was mentioned in Introduction, this theorem earlier was proved under the additional 
assumption that $X$ either has no atoms or has an infinite set of atoms. Such a measure space we 
shall call {\it convenient}. 

We recall the construction of our desired tensor product in the case of a convenient $X$. Take, as 
 underlying linear space of $\Theta$, just $E\ot F$, and as $\theta$ the canonical bilinear 
operator оператор $\vartheta:(x,y)\mt x\ot y$. So, our task is to introduce a suitable norm on 
$L(E\ot F)$. 

We first need an ``extended''  version of our fixed diamond operation, this time between elements 
of amplifications of linear spaces. Namely, for $u\in{L}E, v\in{L}F$ we consider the element $u\di 
v:= \vartheta_\ii(u,v)\in{L}(E\ot F)$. In other words,  this ``diamond operation'' is well defined 
by $\xi x\di\eta y:=(\xi\di\eta)(x\ot y),$ with $\xi,\eta\in L, x\in E, y\in F$.

 An isometry on $L$ will be called {\it proper}, if 
its image is the image of a proper projection. Two isometries will be called {\it disjoint}, if the 
intersection of their images is $\{0\}$.

As is well known (in equivalent terms), if $X$ is convenient, then $L_p(X)$ possesses an infinite 
family of mutually disjoint proper isometries. See, e.g.,~\cite[Cor. 9.12.18]{bog} and 
also~\cite[\S14]{lac} or~\cite[III.A]{woj}. 

The following preparatory statement is crucial in our construction.

\medskip
{\bf Proposition 3.4} (~\cite[Prop. 5.6]{hel}). {\it Let $X$ be convenient. Then every $U\in L(E\ot 
F)$ can be represented as 
\[
a\cd\sum_{k=1}^n I_k\cd(u_k\di v_k),
\]
where $a\in\bb(L),u_k\in LE, v_k\in LF$ and $I_k$ are pairwise disjoint proper isometries on $L$. } 

\medskip
Now we have the right to take $U\in L(E\ot F)$ and assign to it the number

 \[
\|U\|_{pL}:=\inf\left\{\|a\|\left(\sum_{k=1}^n\|u\|^p\|v\|^p\right)^{\frac{1}{p}}\right\},
\]
where the infimum is taken over all possible representations of $U$ in the indicated form. It turns 
out that it is just what we need: 

\medskip
{\bf Theorem 3.5.} (~\cite[Theorem 5.18]{hel}). {\it The function $U\mt\|U\|_{{pL}}$ is a 
${L}$-norm on $E\ot F$, and the pair $(E\ot_{pL}F,\vartheta)$, where $E\ot_{pL}F$ denotes $E\ot F$, 
endowed with the indicated ${L}$-norm, is a $p$-convex tensor product of $E$ and $F$.}  

\medskip
{\bf Remark 3.6}. It was assumed in the cited theorem that $E$ and $F$ are $L$-spaces, and the 
$\di$--operation has the property $\|\xi\di\eta\|=\|\xi\|\|\eta\|$. However one can easily notice 
that its proof uses the estimate $\|a\ot\id_E\|\le\|a\|; a\in\bb(L)$ from Definition 2.1 only for 
operators of rank 1, and only the property of $\di$ to have norm 1. 

\bigskip
We proceed to the main contents of the present paper. How can one behave, if $X$ is not convenient, 
that is the set of its atoms is not empty and finite? It turns out that it is possible to reduce 
the ``unconvenient'' case to the ``convenient'' one. 

For an arbitrary linear space, say $H$, let us consider the algebraic direct sum of a countable 
family of its copies. So, it consists of eventually zero sequences $\bar\xi=(\xi_1,\xi_2,...); 
\xi_k\in H$. If $H$ has a norm, we set $\|\bar\xi\|:=(\sum_k\|\xi_k\|^p)^{\frac{1}{p}}$ and call 
the resulting normed space {\it standard extension of $H$}. 

Now, for our fixed $X$, we denote by ${\Bbb N}X$ the measure space which is the disjoint union of a 
countable family of copies of $X$: ${\Bbb N}X:=X_1\sqcup X_2\sqcup\cdots$. Clearly, {\it ${\Bbb 
N}X$ is convenient.} Therefore the space $L_p({\Bbb N}X)$ satisfies, with ${\Bbb N}X$ in the role 
of $X$, the conditions of Theorem 3.5. 

Let ${\Bbb L}$ be the algebraic direct sum of a countable family of copies of $L$. Then we have the
right to consider on the spaces ${\Bbb L}E$ and ${\Bbb L}F$ the norm of the standard extension of 
$LE$ and $LF$, respectively.

We do not know, whether an arbitrary $L$-space is also a
 ${\Bbb L}$-space with respect to the norm  of the standard extension of the given $L$-norm; may be 
 not. Nevertheless, the following fact is valid. 

\medskip
{\bf Proposition 3.7}. {\it If $E$ is a near-$L$-space, then it is also a near-${\Bbb L}$-space.} 

\smallskip
$\tl$ It is easy to verify that the norm on ${\Bbb L}E$, as well as the norm on $LE$, is a \\
cross--norm with respect to the norm of the underlying space of the given near-$L$-space. 
Therefore, by virtue of Proposition 2.3, it suffices to show that for every ${\bf f}\in{\Bbb L}^*$ 
we have $\|{\bf f}\ot\id_E\|\le\|{\bf f}\|$. In what follows, we omit the easy case $p=1$.

For  $\xi\in L$ and $n\in{\Bbb N}$ we denote by $\bar\xi^n\in{\Bbb L}$ the sequence with the $n$-th 
term $\xi$ and all others zeroes. Introduce the functionals $f_n:L\to\co:\xi\mt{\bf f}(\bar\xi^n)$. 
Fix, for a moment, $n$ and consider an element $\bar u=(u_1,u_2,...)\in{\Bbb L}E$ with $u_n:=\xi x$ 
for some $\xi\in L, x\in E$ and $u_m=0$ for $m\ne n$. We see that ${\bf f}\ot\id_E(\bar u)=\sum_m 
(f_m\ot\id_E)(u_m)$. Since sums of such elements give the whole ${\Bbb L}E$, the same equality is 
valid for all $\bar u\in{\Bbb L}E$.  

But since we know what is $E$, the same Proposition 2.3 gives $\|f_m\ot\id_E\|\le\|f_m\|$, for all 
$m$. Further, it is known (and easy to verify) that $\|{\bf f}\|=(\sum_m\|f_m\|^q)^{\frac{1}{q}}$, 
где $q$ is the number, conjugate to $p$. Therefore for every $\bar v=(v_1,...)\in{\Bbb L}E$ we 
obtain  
\[
 \|({\bf f}\ot\id_E)(\bar v)\|\le\sum_m\|f_m\|\|v_m\|\le
 (\sum_m\|f_m\|^q)^{\frac{1}{q}})(\sum_m\|v_m\|^p)^{\frac{1}{p}})=\|{\bf f}\|\|\bar v\|. \q\tr
 \]

\medskip
The spaces $L$ and ${\Bbb L}$ are connected by the isometry $J:L\to{\Bbb L}:\xi\mt(\xi,0,0,...)$ 
and the coisometry $Q:{\Bbb L}\to L:(\xi,\xi_2,...,\xi_n,...)\mt\xi$; of course, $QJ=\id_L$. For 
every linear space $G$ we shall write $J_G$ instead of $J\ot\id_G:LG\to {\Bbb L}G$ and $Q_G$ 
instead of $Q\ot\id_G:{\Bbb L}G\to LG$. 


Our task is to construct a pair $(\Theta,\theta:E\times F\to\Theta)$, satisfying the conditions of 
Definition 3.1. We shall show that, similarly to Theorem 3.5, {\it we can take $E\ot F$ as the 
underlying linear space of $\Theta$, and the canonical bilinear operator as $\theta$}. 

Where to look for the desired norm on $L(E\ot F)$? 

Using the recipe of Proposition 3.7, we transform $E$ and $F$ into near-${\Bbb L}$-spaces. Then we 
introduce  the bilinear operator $\bar\di:{\Bbb L}\times{\Bbb L}\to{\Bbb L}$ by
\[
\bar\xi\bar\di\bar\eta:=J(Q_E\bar\xi\di Q_F\bar\eta).
\] 
Clearly, it is a diamond operation on ${\Bbb L}$.

But since ${\Bbb N}X$ (though, perhaps, not $X$) is convenient, there exists a $p$-convex tensor 
product of ${E}$ and ${F}$ {\it as that of near-${\Bbb L}$-spaces} with respect to any diamond 
operation on ${\Bbb L}$; in particular, we choose $\bar\di$. Moreover, as the ${\Bbb L}$-space 
$\Theta$ we can take $E\ot F$ with the respective norm on ${\Bbb L}(E\ot F)$) that we shall denote 
by $\|\cd\|_{p{\Bbb L}}$. 

Finally, we introduce a norm on $L(E\ot F)$, induced by the injection $J_{E\ot F}$. In other words, 
for $U\in L(E\ot F)$ we set $\|U\|_{pL}:=\|J_{E\ot F}(U)\|_{p{\Bbb L}}$. 

This norm will turn out to be our desired $L$-norm on the desired tensor product. If there is no 
danger of confusion, we shall omit indices in the notation of the respective norms. 

We must verify the needed requirements.

Take $a\in\bb(L), U\in L(E\ot F)$ and set $\bar a:=JaQ\in\bb({\Bbb L}), \bar U:=J_{E\ot F}(U) 
\in{\Bbb L}(E\ot F)$. Then 
\[
\|a\cd U\|=\|(Ja\ot\id_{E\ot F})(U)\|=\|(JaG)\cd(J_{E\ot F}(U)\|=     
\|\bar a\cd\bar U\|\le\|\bar a\|\|\bar U\|=\|a\|\|U\|,
\]
so that we have the contractibility property. 

If $U_k\in L(E\ot F);k=1,2$ have transversal proper projections $P_k$ in $L$, then ${\Bbb 
P}_k:=JP_kQ$ are transversal proper projections in ${\Bbb L}$, that are supports of $J_{E\ot 
F}U_k$. Therefore the $p$-convexity of $E\ot F$ as an ${\Bbb L}$-space implies the $p$-convexity of 
$E\ot F$ as an ${L}$-space. 

Finally, for $u\in LE, v\in LF$ the equality $J_{E\ot F}(u\di v)=J_Eu\bar\di J_Fv$ and ${\Bbb 
L}$-contractibility of $\vartheta$ imply that $\|u\di v\|\le\|u\|\|v\|$, that is the desired 
$L$-contractibility of $\vartheta$. 

Now the main thing remains: the universal property. In this connection, the following notion will 
be useful. 

\medskip
{\bf Definition 3.8.} Let $G$ be a $p$-convex $L$-space and simultaneously a $p$-convex ${\Bbb 
L}$-space. Then the latter space is called an {\it inflation} of the former space, if $J_G$ is an 
isometry, and $Q_G$ is a coisometry.

\medskip
For example, if $G$ is the minimal $L$-space, then it is easy to show that $G$ as the minimal  
${\Bbb L}$-space is an inflation of the former. As another example, suppose that $G$ belongs to the 
class $SQ_p$, i.e. it is a subspace of a quotient space of some $L_p(Y)$. We make it an $L$-space 
and an ${\Bbb L}$-space by the identification of $LG$ and ${\Bbb L}G$ with the corresponding 
subspaces in $L_p(X,G)$ and $L_p({{\Bbb N}X},G)$. Then the second space is an inflation of the 
first one. The required properties follow from Theorems 1.35 and 1.41 in~\cite{dal}. 

\medskip
In the following three propositions we suppose that $G$ is a given $L$-space that has an inflation,
and we fix the latter.

Let $\rr:E\times F\to G$ be a bilinear operator which is $L$-bounded as an operator between 
near-$L$-spaces with respect to our initial $\di$-operation. We denote by $\bar\rr$ the same 
bilinear operator as an operator between near-${\Bbb L}$-spaces. Speaking about its ${\Bbb 
L}$-boundedness, we mean the diamond operation $\bar\di$, defined above. 

\medskip
{\bf Proposition 3.9.} {\it Our $\bar\rr$ is also ${\Bbb L}$-bounded, and we have 
$\|\bar\rr_\ii\|=\|\rr_\ii\|$.} 

\smallskip
$\tl$ The estimate $\|\bar\rr_\ii\|\le\|\rr_\ii\|$ follows from the formula
\[
 \bar\rr_\ii(\bar u,\bar v)=J_G\Big(\rr_\ii(Q_E\bar u,Q_F\bar v)\Big); \bar u\in{\Bbb L}E, \bar v\in{\Bbb 
L}F, 
\]
which is an easy corollary of the definition of $\bar\di$. The inverse estimate follows from the 
formula 
\[ 
\rr_\ii(u,v)=Q_G\Big(\bar\rr_\ii(J_Eu,J_Fv)\Big),
\]
an easy corollary of the obvious equality $\xi\di\eta=Q(J\xi\bar\di J\eta)$. $\tr$ 

\medskip
Consider the operators $R_\ii:L(E\ot F)\to LG$ and $\bar R_\ii:{\Bbb L}(E\ot F)\to{\Bbb L}G$. These 
are the amplifications of the operator $R:E\ot F\to G$ that is associated with $\rr$ and $\bar\rr$, 
respectively. 

\medskip
{\bf Proposition 3.10.} {\it We have $\|R_\ii\|\le\|\bar R_\ii\|$.} 

\smallskip
$\tl$ This estimate follows from the formula 
\[ 
R_\ii(U)=Q_G\Big(\bar R_\ii(J_{E\ot F}(U))\Big).
\]
Obviously, one should only verify the latter equality on $U$ of the form $\xi(x\ot y); \xi\in 
L,x\in E, y\in F$. Then 
\[
R_\ii(U)=(QJ\xi)\bar R(x\ot y)=Q_G\Big(\bar R_\ii(J\xi(x\ot y)\Big)=Q_G\Big(\bar R_\ii(J_{E\ot F}(U)\Big).\quad \tr 
\]

\medskip
{\bf Proposition 3.11.} {\it If $G$ is as above, then for an arbitrary $L$-bounded bilinear 
operator  $\rr:E\times F\to G$ and the associated linear operator $R:E\ot F\to G$ we have 
$\|R_\ii\|=\|\rr_\ii\|$.} 

\smallskip
$\tl$ Consider $G$ with the ${\Bbb L}$-norm of the given inflation. Because of the universal 
property of the tensor product of our $E$ and $F$ as near-${\Bbb L}$-spaces, we have 
$\|\bar\rr_\ii\|=\|\bar R_\ii\|$. Combining this with two previous propositions, we obtain the 
estimate $\|R_\ii\|\le\|\rr_\ii\|$. 

Further, it follows from the formula $u\di v=Q_{E\ot F}(J_Eu\bar\di J_Fv)$, which can be easily 
verified on elementary tensors, that for all $u\in LE,v\in LF$ we have 
\[
\|\rr_\ii(u,v)\|=\|R_\ii(u\di v)\|\le
\|R_\ii\|\|Q_{E\ot F}(J_Eu\bar\di J_Fv)\|\le\|R_\ii\|\|J_Eu\bar\di J_Fv\|.
\]
But $\vartheta:(x,y)\mt x\ot y$ is ${\Bbb L}$-contractive with respect to the corresponding 
near-${\Bbb L}$-norms and the operation $\bar\di$. Therefore $\|J_Eu\bar\di 
J_Fv\|\le\|J_Eu\|\|J_Fv|\|$. Consequently we have $\|\rr_\ii(u,v)\|\le\|R_\ii\|\|u\|\|v\|$, that is 
$\|\rr_\ii\|\le\|R_\ii\|$. $\tr$ 

\section{Existence of inflations and completion of the proof of the main theorem}

Thus, we see that for concluding the proof of the main theorem it suffices to know that every 
$L$-space has at least one inflation. This for some time we did not know. It is natural to begin 
with the testing of the standard extension of the given $L$-norm. However the existence of 
near-$L$-spaces that are not $L$-spaces (see Remark 2.2) makes one to have doubts; it seems to us 
that it does not fit. 

\bigskip
Nevertheless, by virtue of a recent result of T.~Oikhberg, mentioned in Introduction, one can show 
that inflations do always exist. Indeed it is easily seen that, for arbitrary measure spaces $X$ 
and $Y$ with infinite-dimensional separable $L_p(X)$ and $L_p(Y)$, his argument actually allows us 
 to construct a certain $L_p(Y)$-norm on some $G$, embarking from a given $L_p(X)$-norm on the same $G$.
 We use his method for a proof of the following fact. 

\medskip
{\bf Theorem 4.1}. {\it Let $X$ be as in the formulation of the main theorem. Then for 
$L:=L_p(X);p\in[1,\ii)$ every $L$-space has an inflation.} 

Before the proof, we note that we shall construct an inflation, which essentially differs from the
standard extension of the given $L$-space; see the discussion above.

\smallskip
In what follows, if $Z$ is a measurable subset of some measure space, say $Y$, we shall denote its 
{\it normalized in $L_p(Y)$} characteristic function by $\widehat\chi(Z)$. 
 
\smallskip
$\tl\tl$ We need the following preparatory statement.

\medskip
 {\bf Lemma}. (Here we strictly follow the argument of Oikhberg). {\it Let $Y$ be an 
arbitrary measure space, $\widetilde L$ a subspace in $L_p(Y)$, which is the linear span of several 
characteristic functions of measurable sets. Then there exists  a projection on $\widetilde L$ in 
$\bb((L_p(Y))$ of norm 1.} 

\smallskip
$\tl$ There exist disjoint subsets of non-zero measure $Z_k$ in $Y$, such that $\widetilde 
L=span\{\hat\chi(Z_k)\}$. Take in $L_q(Y)=L_p(Y)^*$ (here $q$ is the conjugate number to $p$) norm 
1 functions  $\widetilde\xi_k$, such that $\widetilde\xi_k=0$ outside $Z_k$, and 
$\la\widetilde\xi_k,\hat\chi(Z_k)\ra=1$. Consider the operator 
$P:L_p(Y)\to L_p(Y):\eta\mt\sum_{k=1}^n\la\widetilde\xi_k,\eta\ra\hat\chi_k$. Clearly, $P$ is 
identical on $\widetilde L$. Further, since for $\eta\in L_p(Y)$ we have 
$P(\eta)=\sum_{k=1}^n\la\widetilde\xi(Z_k),\eta_k\ra\hat\chi(Z_k)$, where $\eta_k=\eta$ on $Z_k$ 
and $\eta_k=0$ outside $Z_k$, we easily see that $\|P(\eta)\|\le\|\eta\|$. $\tr$ 

\bigskip
So, we are given, for $L:=L_p(X)$, an $L$-space $G$. At first we want to introduce a certain norm 
on ${\Bbb L}^0G$, where ${\Bbb L}^0$ is a dense subspace in ${\Bbb L}$, consisting of simple 
functions. 

Every $u\in{\Bbb L}^0G$ can be represented, for some family $Y_k;r=1,...,n$ of pairwise disjoint 
subsets of non-zero measure in ${\Bbb N}X$, as $\sum_{k=1}^n\widehat\chi(Y_k)x_k$,
$x_k\in G$. Take in $X$ an arbitrary family $Z_k$ of pairwise disjoint subsets of non-zero measure 
and set $v:=\sum_{k=1}^n\widehat\chi(Z_k)x_k\in LG$. We put $\|u\|:=\|v\|$. 

\medskip
The subsequent argument consists of several natural stages. 

\medskip
1. {\it The number $\|u\|$ does not depend on the choice of the subsets $Z_k$.} 

\smallskip
$\tl$ Let $Z_k'$ be another family, and $v':=\sum_{k=1}^n\widehat\chi(Z_k')x_k\in LG$.
Consider the operator $J:span\{\widehat\chi(Z_k)\}\to 
span\{\widehat\chi(Z_k')\}:\widehat\chi(Z_k)\mt\widehat\chi(Z_k')$; clearly, it is an isometric 
isomorphism. By the preceding lemma, there exists a projection $P:L\to span\{\widehat\chi(Z_k)\}$ 
of norm 1. Therefore $\|JP\|=1$, so that the contractibility property for $L$-spaces implies 
$\|v'\|=\|JPv\|\le\|v\|$. A similar argument provides the inverse estimate. $\tr$ 

\medskip
2. {\it The number $\|u\|$ does not depend on the representation of $u$ as a sum of elementary 
tensors of the indicated form.} 

\medskip
If we have another representation of our $u$, then, breaking the subsets, corresponding to both 
families, into the same disjoint unions and using the linear independence of the respective 
characteristic functions, we see that both representations lead to the same representation.
To show that the resulting representation gives the same number as the initial one, it suffices, in 
its turn, to show that the number does not change after breaking one of the initial subsets into 
two disjoint subsets of non-zero measure, say, after breaking $Y_1$ into $Y'$ and $Y''$ . Thus, the 
new representation has the form 
$\widehat\chi(Y')z_1+\widehat\chi(Y'')z_2+\sum_{k=2}^n\widehat\chi(Y_k)x_k$ for some $z_1,z_2\in 
G$. Since the tensor factors $\widehat\chi(Y')$ and $\widehat\chi(Y'')$ are linearly independent, 
it follows that $z_l=\lm_lx_1; l=1,2$ for some $\lm_1,\lm_2\in\co$. 

Recall the subsets $Z_k\subset X$. By stage 1, we can assume that we can break $Z_1$ into two 
disjoint sets of non-zero measure, say, $Z'$ and $Z''$. Therefore, if we consider the indicated new 
representation of our $u$, then  the mentioned recipe gives the number $\|v''\|$, where
$v'':=\lm_1\widehat\chi(Z')x_1+\lm_2\widehat\chi(Z'')x_1+\sum_{k=2}^n\widehat\chi(Z_k)x_k$. 

Obviously we have $\|\lm_1\widehat\chi(Z')+\lm_2\widehat\chi(Z'')\|^p=|\lm_1|^p+|\lm_2|^p=1$. This 
easily implies that there exist operators $J_1,J_2$ of norm 1, acting on the space \\
$span\{\widehat\chi(Z'),\widehat\chi(Z''),\widehat\chi(Z_k); k=2,...,n\}$ and such that  
\[
J_1(\widehat\chi(Z_1)=\lm_1\widehat\chi(Z')+\lm_2\widehat\chi(Z'') \q {\rm and} \q 
J_2(\lm_1\widehat\chi(Z')+\lm_2\widehat\chi(Z''))=\widehat\chi(Z_1).
\] 
 Further, the lemma provides a projection $P:L\to 
span\{\widehat\chi(Z'),\widehat\chi(Z''),\widehat\chi(Z_k); k=2,...,n\}$ of norm 1. Therefore, the 
contractibility property for $G$ as for an $L$-space, being applied to the operators $J_1P$ and 
$J_2P$, gives $\|v\|=\|v''\|$. $\tr$ 

\medskip
So, stages 1 and 2 together show that the number $\|u\|; u\in{\Bbb L}^0G$ is well defined.

\medskip
3. {\it The function $u\mt\|u\|$ is an ${\Bbb L}^0$-norm on $G$.} 

\medskip
$\tl$ Obviously, this function is a norm on ${\Bbb L}^0G$. So, it remains to show that for  
$a\in\bb({\Bbb L}^0)$ and $u\in{\Bbb L}^0G$ we have $\|a\cd u\|\le\|a\|\|u\|$. 

Let $u$ has its initial representation, and $a\cd u$ is represented as 
$\sum_{l=1}^m\widehat\chi(Y_l^1)y_l$, where $Y_l^1$ are some pairwise disjoint subsets in ${\Bbb 
N}X$ of non-zero measure. Take $Z_k$ as before and choose arbitrarily one more family 
$Z_l^1;l=1,...,m$ of pairwise disjoint subsets of non-zero measure in $X$. 

Consider the operators $I:span\{\widehat\chi(Y_k)\}\to 
span\{\widehat\chi(Z_k)\}:\widehat\chi(Y_k)\mt\widehat\chi(Z_k)$ and 
$J:span\{\widehat\chi(Y_l^1)\}\to 
span\{\widehat\chi_(Z_l^1)\}:\widehat\chi(Y_l^1)\mt\widehat\chi(Z_l^1)$. Clearly, both are 
isometric isomorphisms. Further, the construction of the norm on ${\Bbb L}^0G$ exactly means that 
$\|u\|=\|(I\ot\id_G)u\|$ and $\|a\cd u\|=\|(J\ot\id_G)(a\cd u)\|$. 

Our lemma provides projections $P:L\to span\{\widehat\chi(Z_k)\}$  and $Q:{\Bbb L}\to 
span\{\widehat\chi(Z_l^1)\}$ of norm 1. We easily see that 
$\Big((P\ot\id_G)(I\ot\id_G)\Big)(u)=\Big(I\ot\id_G\Big)(u)$ and $(Q\ot\id_G)(a\cd u)=a\cd u$.  Now 
set $b:=JQaI^{-1}P:L\to L$. Then we have 
$b\cd\Big((I\ot\id_G)u\Big)=(JQa\ot\id_G)u=(J\ot\id_G)(a\cd u)$. Therefore the contractibility 
property of an $L$-space $G$ gives $\|a\cd 
u\|=\|b\cd\Big((I\ot\id_G)u\Big)\|\le\|b\|\|(I\ot\id_G)u\|\le\|a\|\|u\|$. $\tr$

\medskip
4. {\it There exists an ${\Bbb L}$-norm on $G$ which is an extension of the constructed ${\Bbb 
L}^0$-norm.} 

\smallskip 
$\tl$ Take some $u=\sum_{k=1}^K\xi_kx_k\in{\Bbb L}G; \xi_k\in{\Bbb L},x_k\in G$. 
Consider a sequence $u_n:=\sum_k\xi_{k,n}x_k\in{\Bbb L}^0G$, where $\xi_{k,n}$ converges to $\xi_k$ 
in ${\Bbb L}$. Clearly, the sequence $\|u_n\|$ converges. Denote its limit by $\|u\|$; obviously, 
it does not depend on the choice of $\xi_{k,n}$. 

Take another representation of $u$, say $u=\sum_{l=1}^M\eta_ly_l$, and denote by $\|u\|'$ the 
number, corresponding to this representation. Take linearly independent $z_r\in G; r=1,...,N$, such 
that $x_k=\sum_{r=1}^N\lm_{kr}z_r$ and $y_l=\sum_{r=1}^N\mu_{lr}z_r$ for some 
$\lm_{k,r},\mu_{l,r}\in\co$. We come, in particular, to the representation 
$u=\sum_{r=1}^N(\sum_k\lm_{k,r}\xi_k)z_k$; this leads to some number, denoted by $\|u\|''$.

Recalling $\xi_{k,n}$, we see that $\|u\|''$ is the limit of the sequence
 $\|\sum_{r=1}^N(\sum_k\lm_{k,r}\xi_{k,n})z_r\|$, that clearly coincides with $\|u_n\|$ above. Therefore 
 $\|u\|''=\|u\|$.
 
Similarly, $\|u\|'$ is the number, corresponding to the representation $u$ as 
$\sum_{r=1}^N(\sum_l\mu_{l,r}\eta_l)z_k$. But, since $z_r$ are linearly independent, 
$\sum_k\lm_{k,r}\xi_{k}=\sum_l\mu_{l,r}\eta_l$. Consequently, $\|u\|'$ coincides with $\|u\|''$, 
and hence with $\|u\|$. 

Thus, we have a well defined function $u\mt\|u\|$ on ${\Bbb L}G$, which is obviously a seminorm. 
Let us show that for all $a\in\bb({\Bbb L})$ and $u\in{\Bbb L}G$ we have $\|a\cd u\|\le\|a\|\|u\|$. 

At first suppose that $u\in{\Bbb L}^0G$. Let $u=\sum_{k=1}^N\xi_kx_k$, where $\xi_k\in{\Bbb L}^0$, 
$x_k\in G$. By the lemma, there is a projection of norm 1, say $P$, of ${\Bbb L}$ on a linear span 
of several characteristic functions of measurable sets, such that $u=P\cd u$. 

Take $\epsilon>0$. Since all $a(\xi_k)$ can be approximated by simple functions, the same lemma 
gives a projection of norm 1, say $Q$, on a linear span of several characteristic functions of 
measurable sets, such that $\|a\cd \xi_k-Qa\cd\xi_k\|<\epsilon$. Consider an operator on ${\Bbb 
L}$, acting as $QaP$. Then stage 3 gives $\|(QaP)\cd u\|\le\|QaP\|\|u\|\le\|a\|\|u\|$. Hence 
$\|a\cd u\|\le\|a\cd u-(QaP)\cd u\|+\|(QaP)\cd u\|<\epsilon+\|a\|\|u\|$. This, of course, implies 
$\|a\cd u\|\le\|a\|\|u\|$.

Finally, let $u$ be an arbitrary element in ${\Bbb L}G$, and $u=\sum_{k=1}^N\xi_kx_k$, where 
$\xi_k\in{\Bbb L}$, $x_k\in G$. Since $a$ is bounded, we see that for some sequence $u_n\in{\Bbb 
L}^0G$ we simultaneously have $u_n\to u$ and $a\cd u_n\to a\cd u$; $n\to\ii$. Therefore taking 
limit in the already obtained estimate $\|a\cd u_n\|\le\|a\|\|u_n\|$ we have our desired estimate. 

\medskip
Now let us show that our seminorm is actually a norm.

Take $u\ne0$ and represent it as $\sum_{k=1}^N\xi_kx_k$, with  linearly independent $\xi_k\in{\Bbb 
L}$ and $x_1\ne0$. There exist $a\in\bb({\Bbb L})$ such that $a(\xi_1)=\xi_1$ and $a(\xi_k)=0; 
k>1$, and also, for every $k=1,...,N$ a sequence $\xi_{k,n}\in{\Bbb L}^0; n\in{\Bbb N}$, converging 
to $\xi_k$. Set $u_n:=\sum_{k=1}^N\xi_{k,n}x_k$; we have $\|a\cd 
u_n\|\ge\|a(\xi_{1,n})x_1\|-\sum_{k=2}^N\|a(\xi_{k,n})x_k\|$. But clearly $a(\xi_{1,n})$ converges 
to $\xi_1$, and $a(\xi_{k,n})$ converges to 0 for other $k$. Hence for sufficiently big $n$ we have 
$\|a\cd u_n\|\ge\epsilon$ for some $\epsilon>0$. Combining this with the estimation above, we 
obtain that $\|u_n\|\ge\epsilon/\|a\|$.Therefore $\|u\|>0$. $\tr$ 

\medskip
5. {\it The $p$-convexity is preserved by passing from $L$-spaces to ${\Bbb L}$-spaces.}

\smallskip 
$\tl$ Suppose, at first, that $u_1,u_2$ with transversal supports lie in ${\Bbb L^0}G$. By 
definition of the norm in ${\Bbb L^0}G$, there exist a family $Z_1^1,...,Z_N^1,Z_1^2,...,Z_M^2$ of 
pairwise disjoin subsets of non-zero measure in $X$ such that, for some $v_1,v_2\in LG$ of the form 
$v_1:=\sum_{k=1}^N\widehat\chi(Z_k^1)x_k$, $v_2:= \sum_{k=1}^M\widehat\chi(Z_l^2)y_l$, 
respectively, we have $\|u_1\|=\|v_1\|$, $\|u_2\|=\|v_2\|$ and $\|u_1+u_2\|=\|v_1+v_2\|$. Since 
$v_1,v_2$ obviously have transversal supports, and our $L$-space is $p$-convex, we have 
$\|u_1+u_2\|^p\le\|u_1\|^p+\|u_2\|^p$. 

Now take arbitrary $u_1,u_2\in{\Bbb L}G$ with transversal supports. Clearly for $k=1,2$  there 
exist a sequence $u_k^n\in{\Bbb L^0}G$ with the same support as $u_k$, converging to $u_k$. Then, 
passing to limits, we obtain the desired estimate. $\tr$ 

\bigskip 
It is clear that all statements and arguments in stages 1-5 are valid, if we replace ${\Bbb 
L}:=L_p({\Bbb N}X)$ by $L_p(Y)$ for an arbitrary measure space $Y$. Now we concentrate on our 
concrete situation.

\medskip
{\bf End of the proof.} It remains to show that $J_G$ is an isometry, and $Q_G$ is a coisometry.

$\tl$ Take at first $u\in L^0G$ and represent it as $\sum_{k=1}^N\widehat\chi(Y_k)x_k$ with 
pairwise disjoint $Y_k\subset X$. Then $J_G(u)$, as an element of the subspace ${\Bbb L}^0G$ of 
${\Bbb L}G$ has the same representation, only now we must consider $Y_k$ as subsets in the first 
summand $X_1$ in ${\Bbb N}X= X_1\sqcup X_2\sqcup...$. Therefore, calculating $\|J_G(u)\|$ by the 
prescribed recipe, we can take as $Z_k$ the initial $Y_k$, and the same $x_k$. But then the 
respective $v$ is just $u$, therefore $\|J_G(u)\|=\|u\|$. Thus, the restriction of $J_G:LG\to{\Bbb 
L}G$ on a dense subspace in $LG$ is an isometry, so the same is true for $J_G$. 

Turn to $Q_G$. Since we have $Q_GJ_Gu=u$ for all $u\in LG$, and $J_G$ is an isometry, we only have 
to show that $Q_G$ 
is contractive. 

Take $\bar u\in{\Bbb L}$ and observe that $J_GQ_G\bar u=(JQ)\cd\bar u$. Therefore $\|Q_G\bar 
u\|=\|J_GQ_G\bar u\|\le\|JQ\|\|\bar u\|=\|\bar u\|$. $\tr$ 

\bigskip
The proof of Theorem 4.1 is concluded. $\tr\tr$

\bigskip
Combining this theorem with Proposition 3.11, we obtain our main theorem.

\section {$p$-convex tensor product as a functor} 

\bigskip
Now let us do some final observation. Recall an important notion in geometry of normed spaces. 
Suppose that we assign to every pair, say $E,F$, of normed spaces the space $E\ot F$ endowed with 
some norm. The most interesting are assignments, satisfying the so-called {\it metric mapping 
property}~\cite[\S 12]{def} (see also ``uniform cross-norms'' in~\cite[\S 6.1]{rya}): for every 
bounded operators $\va:E_1\to F_1, \psi:E_2\to F_2$ we have 
 $\|\va\ot\psi\|\le\|\va\|\|\psi\|$. (In other terms, such an assignment, extended to bounded 
 operators, is a bifunctor on the category of normed spaces and contractive operators.) We shall show that the $p$-convex tensor 
 product has a natural analogue of that ``functorial property'' for near-$L$-spaces. As usual, 
 $L:=L_p(X)$, and the only condition on $X$ is that $L$ is infinite-dimensional and 
 separable.

\medskip
{\bf Proposition 5.1}. {\it Let $\va:E_1\to E_2$, $ \psi:F_1\to F_2$ be $L$-bounded operators 
between near-$L$--spaces. Then the operator $\va\ot\psi:E_1\ot_{pL}E_2\to F_1\ot_{pL}F_2$ is 
$L$--bounded, and  we have $\|(\va\ot\psi)_\ii\|\le\|\va_\ii\|\|\psi_\ii\|$.} 

\smallskip
$\tl$ We first assume that $X$ is convenient. Then, taking $U\in L(E_1\ot E_2)$, we have the right 
to present it as $a\cd\sum_k I_k\cd u_k\di v_k; a\in\bb(L), u_k\in LE_1,v_k\in LF_1$ with mutually 
disjoint proper isometries $I_k$ (see Proposition 3.4). Then, since amplifications of our operators 
are morphisms of left $\bb(L)$-modules, we have 
$(\va\ot\psi)_\ii(U)=a\cd\sum_kI_k\cd(\va\ot\psi)_\ii(u_k\di v_k)$. By the formula 
$(\va\ot\psi)_\ii(u\di v)=\va_\ii(u)\di\psi_\ii(v)$, easily verified on elementary tensors, the 
latter expression is $a\cd\sum_kI_k\cd(\va_\ii(u_k)\di\psi_\ii(v_k)$. Therefore, by definition of 
the norm on $F_1\ot_{pL}F_2$ for convenient $X$, we have  
\[
\|(\va\ot\psi)_\ii(U)\|\le\|a\|\Big(\sum_k\|\va_\ii(u_k)\|^p\|\psi_\ii(v_k)\|^p\Big)^{\frac{1}{p}}\le
\]
\[ 
\|a\|\Big(\sum_k(\|\va_\ii\|\|u_k\|)^p\|(\|\psi_\ii\|\|v_k\|)^p\Big)^{\frac{1}{p}}\le 
\|a\|\|\va_\ii\|\|\psi_\ii\|\Big(\sum_k\|u_k\|)^p\|v_k\|^p\Big)^{\frac{1}{p}}. 
\]

It remains to take the respective infimum over all representations of $U$ in the prescribed form. 

 Turn to an arbitrary $X$. Using the standard extension of given $L$-norms, consider 
 our four spaces as near-${\Bbb L}$-spaces (see Proposition 3.7). Denote our given operators as acting 
 between near-${\Bbb L}$-spaces as $\widetilde\va$ and $\widetilde\psi$, respectively. Thus, since 
 ${\Bbb L}=L_p({\Bbb N}X)$, and ${\Bbb N}X$ is convenient, we have
$\|(\widetilde\va\ot\widetilde\psi)_\ii\|\le\|\widetilde\va_\ii\|\|\widetilde\psi_\ii\|$. But, 
using the definition of near-${\Bbb L}$-norms on our four spaces, we easily obtain that 
$\|\widetilde\va_\ii\|=\|\va_\ii\|$ and $\|\widetilde\psi_\ii\|=\|\psi_\ii\|$. At the same time for 
$U\in L(E_1\ot_{pl}E_2)$ we have 
\[
\|(\va\ot\psi)_\ii(U)\|=\|J_{F_1\ot F_2}\Big(\va\ot\psi)_\ii(U)\Big)\|=
\|\Big((\va\ot\psi)_\ii(U),0,0,...\Big)\|=
\]
\[
\|(\widetilde\va\ot\widetilde\psi)_\ii(U,0,0,...)\|\le\|(\widetilde\va\ot\widetilde\psi)_\ii\|\|(U,0,0,...)\|=
\|(\widetilde\va\ot\widetilde\psi)_\ii\|\|U\|,
\]
and consequently $\|(\va\ot\psi)_\ii\|\le\|(\widetilde\va\ot\widetilde\psi)_\ii\|$. The desired 
estimate immediately follows. $\tr$

A.~Ya.~Helemskii

Faculty of Mechanics and Mathematics

Moscow State (Lomonosov) University

Moscow 119991 Leninskie Gory

E-mail: helemskii@rambler.ru

\end{document}